\documentclass[12pt]{amsart}
\usepackage{amssymb, psfig, epsf}

\theoremstyle{plain}
\newtheorem{thm}{Theorem}
\newtheorem{prop}[thm]{Proposition}
\newtheorem{lemma}[thm]{Lemma}
\newtheorem{cor}[thm]{Corollary}

\theoremstyle{definition}
\newtheorem{exa}[thm]{Example}
\newtheorem{rem}[thm]{Remark}

\newcommand{\bZ}{\mbox{${\mathbb Z}$}}

\newcommand{\frakS}{{\mathfrak S}}

\newcommand{\lI}{\begin{picture}(3,3)
   \put(0,0){\line(1,0){3}}  \put(0,0){\line(0,1){3}}
     \put(0,3){\line(1,0){3}}  \put(3,0){\line(0,1){3}}
 \end{picture}}

\newcommand{\lII}{\begin{picture}(6,6)(0,1)
  \put(0,0){\line(1,0){3}} \put(0,3){\line(1,0){6}} \put(0,6){\line(1,0){6}}
  \put(0,0){\line(0,1){6}} \put(3,0){\line(0,1){6}} \put(6,3){\line(0,1){3}}
 \end{picture}}

\newcommand{\lIII}{\begin{picture}(6,6)(0,1)
   \put(0,0){\line(1,0){3}}  \put(0,0){\line(0,1){3}}
   \put(0,3){\line(1,0){6}}  \put(3,0){\line(0,1){6}}
   \put(3,6){\line(1,0){3}}  \put(6,3){\line(0,1){3}}
 \end{picture}}

\newcommand{\lIV}{\begin{picture}(3,6)
  \put(0,0){\line(1,0){3}}  \put(0,0){\line(0,1){6}}
  \put(0,3){\line(1,0){3}}  \put(3,0){\line(0,1){6}}
  \put(0,6){\line(1,0){3}}
 \end{picture}}

\newcommand{\lV}{\begin{picture}(6,3)
  \put(0,0){\line(1,0){6}}  \put(0,0){\line(0,1){3}}
  \put(0,3){\line(1,0){6}}  \put(3,0){\line(0,1){3}}
  \put(6,0){\line(0,1){3}}
 \end{picture}}

\newcommand{\spa}{\;\:}

\newcommand{\fs}{{\mathfrak S}}

\newcounter{FNC}[page]
\def\newfootnote#1{{\addtocounter{FNC}{2}$^\fnsymbol{FNC}$%
     \let\thefootnote\relax\footnotetext{$^\fnsymbol{FNC}$#1}}}

\begin{document}
\bibliographystyle{amsplain}
\title[Skew Schubert polynomials]{Skew Schubert polynomials}

\author{Cristian Lenart}
\address{Department of Mathematics and Statistics\\
         State University of New York at Albany\\
         Albany, NY \ 12222\\
         USA}
\email{lenart@csc.albany.edu}

\author{Frank Sottile}
\address{Department of Mathematics\\
        University of Massachusetts\\
        Amherst, MA \ 01003\\
        USA}
\email{sottile@math.umass.edu}
\date{5 May 2002}
\thanks{Most of this work was done while the first author was supported by Max-Planck-Institut f\"ur Mathematik. Second author supported in part by NSF grants  DMS-9701755 and
         DMS-0070494.}

\begin{abstract}
We define skew Schubert polynomials to be normal form (polynomial)
representatives of certain classes in the cohomology of a flag manifold.
We show that this definition extends a recent construction of Schubert
polynomials due to Bergeron and Sottile in terms of certain increasing
labeled chains in Bruhat order of the symmetric group.
These skew Schubert polynomials
expand in the basis of Schubert polynomials with nonnegative integer
coefficients that are precisely the structure constants of the cohomology of
the complex flag variety with respect to its basis of Schubert classes.
We rederive the construction of Bergeron and Sottile in a purely
combinatorial way, relating it to the  construction of Schubert
polynomials in terms of rc-graphs.
\end{abstract}

\maketitle

\section*{Introduction}

Skew Schur polynomials $S_{\lambda/\mu}(x_1,\dotsc,x_k)$ play an
important and multi-faceted role in algebraic combinatorics.
For example, they are generating functions for the Littlewood-Richardson
coefficients $c^\lambda_{\mu,\nu}$
 \begin{equation}\label{E:LR-Schur}
   S_{\lambda/\mu}(x_1,\dotsc,x_k)\ =\
    \sum_\nu  c^\lambda_{\mu,\nu}\, S_\nu(x_1,\dotsc,x_k)\,.
 \end{equation}
They are also generating functions for Young tableaux of shape
$\lambda/\mu$
 \begin{equation}\label{E:Comb-Schur}
   S_{\lambda/\mu}(x_1,\dotsc,x_k)\ =\
    \sum  {\bf x}^T\,,
 \end{equation}
where the sum is over all Young tableaux $T$ of shape $\lambda/\mu$ and
${\bf x}^T$ is a monomial associated to $T$.
The relationship between these two very different facets of skew Schur
polynomials involves the combinatorics of Young
tableaux~\cite[Chapter 7, Appendix 1]{St99}.
Young tableaux can be thought of
as increasing labeled chains in Young's lattice, with covers labeled by pairs
$(k,l)$,
where $k$ is the entry in the box corresponding to the cover, and
$l$ the content of that box (that is, the difference between its
column and row); labels are ordered lexicographically.
We introduce skew Schubert polynomials,
which are the Schubert polynomial analogs of skew Schur
polynomials, in the sense that they
generalize the two properties mentioned above.


In Section 1, we introduce our key concepts of increasing labeled chains in the Bruhat
order, which generalize Young tableaux, and of skew Schubert polynomials, $\frakS_{w/u}(x)$
for permutations $u\leq w$.
In particular, we show that skew Schubert polynomials are generating functions
both for Littlewood-Richardson coefficients $c^w_{u,v}$ for Schubert
polynomials (the analog of~\eqref{E:LR-Schur}), and for increasing chains in
the Bruhat order (the analog of~\eqref{E:Comb-Schur}).
Our skew Schubert polynomials
are different from those defined in \cite{Ki97} (using skew
divided difference operators), which do not expand with
nonnegative coefficients in the basis of Schubert polynomials, in general.

The skew Schubert polynomial $\frakS_{w_0/w}(x)$ is the ordinary Schubert
polynomial $\frakS_w(x)$, and our formula for $\frakS_{w_0/w}(x)$ in terms of
increasing chains is the formula of Corollary~5.3 in~\cite{BS02}, which
inspired us (see the next section for the notation).
In Section 2, we give a purely combinatorial proof of that
formula, 
relating it to a standard
construction~\cite{FK96, BJS93, FS95} in terms of rc-graphs by giving a
content-preserving bijection between rc-graphs and increasing chains.

\section{Skew Schubert Polynomials}

Let $\Sigma_n$ be the symmetric group of permutations of $\{1,2,\dotsc,n\}$.
For $w\in\Sigma_n$, the \emph{length}, $\ell(w)$, of $w$ is the number of
inversions of $w$.
Set $w_0\in\Sigma_n$ to be the longest permutation,
$w_0:=n\,\dotsc\,2\,1$.
We use basics on Schubert and Schur polynomials, which may be found in any
of~\cite{LS82a,Mac91,Fu97,Man98}.

The main outstanding problem in the theory of Schubert polynomials is the
Littlewood-Richardson problem~\cite[Problem 11]{St01}:
Determine the structure constants
$c^w_{u,v}$ defined by the polynomial identity
 \[
    \frakS_u(x) \cdot \frakS_v(x)\ =\ \sum_w c^w_{u,v} \frakS_w(x)\,.
 \]
Since every Schur polynomial is a Schubert polynomial, this problems asks for
the analog of the classical Littlewood-Richardson rule.
The classical Littlewood-Richardson coefficients $c^\lambda_{\mu,\nu}$ have
important and intricate combinatorial
properties~\cite[Chapter 7, Appendix 1]{St99}.
The Littlewood-Richardson coefficients $c^w_{u,v}$ for
Schubert polynomials should be similarly important.
Indeed, they are intersection numbers of Schubert varieties; more
precisely, $c^w_{u,v}$ enumerates flags in a suitable triple
intersection of  Schubert varieties indexed by $u$, $v$, and $w_0w$.

The cohomology classes $\sigma_w$ (called Schubert classes) for $w\in\Sigma_n$
of Schubert varieties form an integral basis for the cohomology ring
$H^*{\mathbb F}\ell_n$ of the flag manifold.
Many identities and some formulas for the constants $c^w_{u,v}$ have been
obtained by studying the class $\sigma_u\cdot\sigma_{w_0w}$
in $H^*{\mathbb F}\ell_n$~\cite{So96,BS98,BS_lag-pieri}.
This is because, for $u,v,w\in\Sigma_n$, we have the following identity in
the cohomology ring
 \[
   \sigma_u\cdot\sigma_{w_0w}\cdot\sigma_v\ =\ c^w_{u,v}\,\sigma_{w_0}
 \]
and hence, by the duality of the intersection pairing
 \begin{equation}\label{coh-formula}
   \sigma_u\cdot\sigma_{w_0w}\ =\
      \sum_{v\in\Sigma_n} c^w_{u,v}\,\sigma_{w_0v}\,.
 \end{equation}

With this motivation, we would like to define the skew Schubert polynomial
$\frakS_{w/u}(x)$ for permutations $u\leq w$ in $\Sigma_n$
to be \emph{the} polynomial representative of the class
$\sigma_u\cdot\sigma_{w_0w}$.
Unfortunately, the cohomology of the flag manifold is
isomorphic to the quotient ring
 \[
    H^*{\mathbb F}\ell_n\ :=\ {\mathbb Z}[x_1,x_2,\ldots,x_n]/
    \langle e_i(x_1,\ldots,x_n)\mid i=1,\ldots,n\rangle\,,
 \]
and so there is no well-defined polynomial representative of
$\sigma_u\cdot\sigma_{w_0w}$.
Here $e_i(x_1,\dotsc,x_n)$ is the $i$th elementary symmetric polynomial in
$x_1,\dotsc,x_n$.

If however we choose the degree reverse lexicographic term order on
monomials in the polynomial ring ${\mathbb Z}[x_1,x_2,\ldots,x_n]$, with
$x_1<x_2<\dotsb<x_n$, then every
element of this quotient ring has a unique normal form
polynomial representative~\cite[Prop.~1.1]{Sturmfels_GBCP} with respect to
this term order.
These normal form representatives are elements of
 \[
   \mathbb{Z}\cdot\{x^\alpha\mid \alpha=(\alpha_1,\dotsc,\alpha_n)\
     \text{ with }\ \alpha_i\leq n-i\}\,.
 \]
Equivalently, if we set $\delta:=(n-1,\dotsc,2,1,0)$, then normal form
representatives of cohomology classes are sums of monomials dividing
$x^\delta$.

Fomin, Gelfand, and Postnikov~\cite{FGP97} observed that for $w\in\Sigma_n$,
the Schubert polynomial $\frakS_w(x)$ is the normal form
representative of the corresponding Schubert class $\sigma_w$.
In fact, this feature of Schubert polynomials, that one does not need to work
modulo an ideal, is what led Lascoux and Sch\"utzenberger to their definition
of Schubert polynomial.
Thus we define the {\it skew Schubert polynomial} $\frakS_{w/u}(x)$ to
be the normal form representative
of the class $\sigma_u\cdot\sigma_{w_0w}$.
These skew Schubert polynomials are generating
functions for the coefficients $c^w_{u,v}$.

\begin{thm}\label{TheoremOne}
 For $u\leq w$ in $\Sigma_n$, we have
 \begin{equation}\label{E:thm1}
   \frakS_{w/u}(x)\ =\ \sum_v c^w_{u,v}\,\frakS_{w_0v}(x)\,.
 \end{equation}
\end{thm}

\begin{proof}
 Each side of~\eqref{E:thm1} is the normal form representative of the
 corresponding side of~\eqref{coh-formula}.
\end{proof}

While $\frakS_{w/u}(x)$ depends upon $n$, the Laurent polynomial
$x^{-\delta}\frakS_{w/u}(x)$ does not, as
$\frakS_u(x)$ and $x^{-\delta}\frakS_{w_0w}(x)$ are each independent of $n$.
The independence of $\frakS_u(x)$ is the familiar stability property of
Schubert polynomials.
The independence of $x^{-\delta}\frakS_{w_0w}(x)$
may be similarly deduced using divided differences.
This also follows from Proposition~\ref{prop:Schub-constr}, which we prove
in Section 2.
\smallskip

The skew Schubert polynomial $\frakS_{w/u}(x)$ is also a generating function for
certain chains from $u$ to $w$ in the Bruhat order.
The Bruhat order is defined by its covers; $u\lessdot w$ if and only if
$u^{-1}w$ is a transposition $(k,l)$ with $u(k)<u(l)$, and for every
$k<i<l$ we have either $u(i)<u(k)$ or $u(l)<u(i)$.
Thus $u\lessdot u\cdot(k,l)$ with $k<l$ if $u(k)<u(l)$ and at no position
between $k$ and $l$ does $u$ take a value between $u(k)$ and $u(l)$.
This implies that the difference in lengths $\ell(w)-\ell(u)$
equals 1.

The {\it labeled Bruhat order} is the labeled r\'eseau (labeled directed
multigraph) obtained from the Bruhat order by drawing a directed edge
$u\xrightarrow{\,(k,b)\,} w$ whenever $u\lessdot w$ with
$u^{-1}w=(i,j)$, where $i\leq k< j$ and $b=u(i)=w(j)$.
We note that all the constructions
and statements still hold if we uniformly set $b=u(j)=w(i)$,
but it is more convenient to use the first definition.
There will be $j{-}i$ such directed edges for every cover.
This structure was defined in~\cite{BeBi} and has been crucial in
subsequent work on the problem of multiplying Schubert polynomials by
Bergeron and Sottile.
In~\cite{BS98} this structure was called the colored Bruhat order.

To any (saturated) chain $\gamma$ in this r\'eseau, we associate a monomial
$x^\gamma$ in the variables $x_1,\ldots,x_{n-1}$, where the power of $x_i$
counts how often $i$ was the first coordinate of a label in $\gamma$.
A chain
 \begin{equation}\label{labchain}
   u_0\ \xrightarrow{\,(k_1,b_1)\,}\
   u_1\ \xrightarrow{\,(k_2,b_2)\,}\ \cdots\
   \xrightarrow{\,(k_m,b_m)\,}\ u_m
\end{equation}
in this r\'eseau is {\it increasing} if its sequence of labels is increasing
in the lexicographic order on pairs of integers.

These definitions allow a nice reformulation of a
combinatorial construction of Schubert polynomials given in~\cite{BS02}.

\begin{thm}\label{thm:skew}
 Let $u\leq w$ be permutations in $\Sigma_n$.
 Then
 \begin{equation}\label{E:Comb-Schubert}
   \frakS_{w/u}(x)\ =\ \sum x^\delta/x^\gamma,
 \end{equation}
 the sum over all increasing chains $\gamma$ in the labeled Bruhat order from
 $u$ to $w$.
\end{thm}

Note that $\frakS_w(x)=\frakS_{w_0/w}(x)$.
We recover the formula of Bergeron and Sottile for $\frakS_w(x)$~\cite{BS02}.
Denote by ${\Gamma}(w,w_0)$ the set of
all increasing labeled chains from $w$ to $w_0$.

\begin{prop}[Bergeron-Sottile~\cite{BS02}]\label{prop:Schub-constr}
 Let $w\in \Sigma_n$. Then
 $$
   \frakS_w(x)\ =\ \sum_{\gamma\in{\Gamma}(w,w_0)} x^\delta/x^\gamma.
 $$
\end{prop}

Theorem~\ref{thm:skew} has an immediate enumerative consequence.
The {\it type} $\alpha$ of a chain $\gamma$ in the labeled Bruhat order
is the (weak) composition $\alpha$ whose $i$th component counts the number of
occurrences of $i$ as the first coordinate of an edge label.
Thus $x^\gamma=x^\alpha$.
For $u\leq w$ in the Bruhat order and  a composition $\alpha$ of
$\ell(w)-\ell(u)$, let $I_\alpha(u,w)$ count the number of
increasing chains from $u$ to $w$ of type $\alpha$.
Since $\frakS_{w_0v}(x)=\frakS_{w_0/w_0v}(x)$, equating coefficients
in~\eqref{E:thm1} and using Theorem~\ref{thm:skew} we obtain

\begin{cor}\label{cor:identity}
  \qquad ${\displaystyle I_\alpha(u,w)\ =\
           \sum_v c^w_{u,v} I_\alpha(w_0v,w_0)}$.
\end{cor}

For a composition $\alpha$ with $n{-}1$ parts where $\alpha\leq \delta$
coordinatewise, set
 \[
   h_\alpha(x)\ :=\ h_{\alpha_1}(x_1)h_{\alpha_2}(x_1,x_2)\dotsm
                    h_{\alpha_{n-1}}(x_1,\ldots,x_{n-1})\,,
 \]
where $h_a(x_1,\ldots,x_k)$ is the complete homogeneous symmetric polynomial
of degree $a$ in the variables $x_1,\ldots,x_k$.
We let $h_a$ and $h_\alpha$ denote the corresponding cohomology classes in
$H^*{\mathbb F}\ell_n$.
For $f\in H^*{\mathbb F}\ell_n$, set $\psi_\alpha(f)$ to be the coefficient of
$\sigma_{w_0}$ in the product $f\cdot h_\alpha$.
This gives a linear map
$\psi_\alpha\colon H^*{\mathbb F}\ell_n\rightarrow{\mathbb Z}$.

\begin{lemma}\label{lem:linear}
 For any $f\in  H^*{\mathbb F}\ell_n$,
 $\psi_\alpha(f)$ is the coefficient of the monomial
 $x^\delta/x^\alpha$ in the normal form representative of $f$.
\end{lemma}

\begin{proof}
 This follows from two observations.
 First, both $\psi_\alpha$ and the
 map associating to $f$ the coefficient of the monomial
 $x^\delta/x^\alpha$ in the normal form representative of $f$
 are ${\mathbb Z}$-linear maps on $H^*{\mathbb F}\ell_n$.
 Second, the lemma holds on the basis of Schubert classes,
 whose representatives are Schubert polynomials.
 This was shown by Kirillov and Maeno~\cite{KM00} and is also a
 consequence of the Pieri formula~\cite{So96} and the construction of Schubert
 polynomials given in Proposition~\ref{prop:Schub-constr}.
\end{proof}

\begin{proof}[Proof of Theorem~\ref{thm:skew}]
 Given a chain $\gamma$ in the Bruhat order, let end$(\gamma)$ denote its end
 point.
 The Pieri-type formula for Schubert classes~\cite{So96} is
 \begin{equation}\label{E:Pieri}
   \sigma_u\cdot h_a\ =\
   \sum \sigma_{\mathrm{end}(\gamma)}(x)\,,
 \end{equation}
 the sum over all increasing chains that begin at $u$, have length $a$, and
 whose labels have first coordinate $k$.

 Applying Lemma~\ref{lem:linear} to the class
 $\sigma_u\cdot\sigma_{w_0w}$ shows that the coefficient
 of the monomial $x^\delta/x^\gamma$ in the normal form representative of
 $\sigma_{w/u}$ is the coefficient of
 $\sigma_{w_0}$ in the triple product
 $\sigma_u\cdot\sigma_{w_0w}\cdot h_\alpha$, where $\alpha$ is the type of the
 chain $\gamma$.
 Since the coefficient of $\sigma_{w_0}$ in a product
 $\sigma_v\cdot\sigma_{w_0w}$ is the Kronecker delta
 $\delta^v_w$ (by the duality of the intersection pairing), this
 is the coefficient of $\sigma_w$ in the product
 $\sigma_u\cdot h_\alpha$.
 The theorem now follows by expanding this product and iteratively applying
 the Pieri formula.
\end{proof}

\begin{rem}
 By the definition of skew Schubert polynomials, we have
 \[
    \frakS_{w_0/w_0w}(x)\ =\ \frakS_{w_0w}(x)\ =\ \frakS_{w/1}(x)\,,
 \]
 since $w_0w_0=1$ and $\frakS_1(x)=1$.
 Thus there should exist a natural bijection between $\Gamma(1,w)$ and
 $\Gamma(w_0w,w_0)$ that preserves type.
 Corollary~\ref{cor:identity} suggests that, for $u\leq w$, this bijection
 should generalize to a map from $\Gamma(u,w)$ to the
 union $\coprod_v\Gamma(w_0v,w_0)$ such that the cardinality of the inverse
 image of each chain in $\Gamma(w_0v,w_0)$ is $c_{u,v}^w$.
 This would give a combinatorial interpretation for
 $c_{u,v}^w$, and thus solve the Littlewood-Richardson problem.
\end{rem}

\begin{rem}
 The skew Schubert polynomial $\frakS_{w/u}(x)$ also has a geometric
 interpretation as the representative in cohomology of a skew Schubert variety
 \[
    X_u F\!\!_{\mbox{\Huge .}}\,\cap\, X_{w_0w}F'\!\!\!_{\mbox{\Huge .}}\,,
 \]
 where $X_w F\!\!_{\mbox{\Huge .}}$ is the Schubert variety whose representative
 in cohomology is $\frakS_w(x)$ and $F\!\!_{\mbox{\Huge .}}$ and
 $F'\!\!\!_{\mbox{\Huge .}}$ are flags in general position.
 Stanley noted that skew Schur polynomials have a similar
 geometric-cohomological interpretation~\cite[\S 3]{St77}.
\end{rem}

\begin{rem}
 Schubert polynomials $\frakS_w(x)$ when $w$ has a single descent at position
 $k$ ($w$ is a Grassmannian permutation)
 are Schur symmetric polynomials in $x_1,\dotsc,x_k$.
 Even if $w$ and $u$ are Grassmannian permutations with the
 same descent, then it is not necessarily the case that the skew Schubert
 polynomial $\frakS_{w/u}(x)$ equals the corresponding skew Schur polynomial.
 We illustrate this in $\Sigma_4$, using partitions and skew partitions for
 indices of Schur polynomials.

 We have the Schubert/Schur polynomials
 \[
   \begin{array}{rclcl}
    \frakS_{1324}(x)&\ =\ &x_1+x_2           &\ =\ &S_{\lI}(x_1,x_2)\,,\quad
           \text{\ and}\\
    \frakS_{2413}(x)&\ =\ &x_1^2x_2+x_1x_2^2 &\ =\ &S_{\lII}(x_1,x_2)\,,
        \rule{0pt}{14pt}
   \end{array}
 \]
 and also the skew Schur polynomial
 \[
   \begin{array}{rcccccl}
    S_{\lIII}(x_1,x_2)&\ =\ &x_1x_2 &+& x_1x_2+x_1^2+x_2^2\\
       &\ =\ & S_{\lIV}(x_1,x_2) &+& S_{\lV}(x_1,x_2)
         &\ =\ &\frakS_{2314}(x) + \frakS_{1423}(x)\,.\rule{0pt}{14pt}
   \end{array}
 \]
 On the other hand,
 $\frakS_{2413/1324}(x)=\frakS_{1324}(x)\cdot\frakS_{w_02413}(x)$
 and $w_02413=3142$.
 By Monk's formula (or the Pieri formula~\eqref{E:Pieri}), we compute
 $\frakS_{2413/1324}(x)$ to be
\newcommand{\hsm}{\hspace{4.7pt}}
 \[
   \begin{array}{rcl}
    \frakS_{1324}(x)\cdot\frakS_{3142}(x)&\ =\ &
      \hspace{5pt}\frakS_{3241}(x)\hsm
          +\hsm\frakS_{4132}(x)\hsm+\hsm\frakS_{3412}(x)\\
      &\ =\ &
        \frakS_{w_02314}(x)+\frakS_{w_01423}(x)+\frakS_{w_02143}(x)\,.
         \rule{0pt}{14pt}
   \end{array}
 \]
 The first two terms of $\frakS_{2413/1324}(x)$ carry the same information as the
 two terms of $S_{\lIII}(x_1,x_2)$, but there is also a third term.
 This is because the multiplication of Schubert polynomials is
 richer than the multiplication of Schur polynomials.
 By Theorem~\ref{TheoremOne}, this expansion of $\frakS_{2413/1324}(x)$
 records the fact that $\frakS_{2413}(x)$ is a summand (with
 coefficient 1) in each of
 \[
   \frakS_{1324}(x)\cdot\frakS_{2314}(x)\,,\qquad
   \frakS_{1324}(x)\cdot\frakS_{1423}(x)\,,\qquad\text{and}\qquad
   \frakS_{1324}(x)\cdot\frakS_{2143}(x)\,.
 \]
 Note that the last product (unlike the first two) is not
 symmetric in $x_1$ and $x_2$ even though both $\frakS_{1324}(x)$
 and $\frakS_{2413}(x)$ are symmetric in $x_1$ and $x_2$.
\end{rem}

\section{RC-Graphs and Increasing Labeled Chains}

RC-graphs are combinatorial objects associated to permutations $w$ in the
symmetric group $\Sigma_n$.
They were defined in \cite{FK96, BeBi} as certain subsets of
 \[
     \{(k,b)\in\bZ_{>0}\times\bZ_{>0}\mid k+b\le n\}\,.
  \]
We may linearly order such a subset of pairs by setting
 \[
   (k,b)\le (j,a) \;\;\Longleftrightarrow\;\;  (k< j) \spa\mbox{or} \spa
   (k=j\spa\mbox{and}\spa b\ge a)\,.
 \]
Let $(k_i,b_i)$ be the $i$th pair in this linear order.
Consider the sequence
 \[
   d(R)\ :=\ (k_1+b_1-1, k_2+b_2-1,\ldots)\,.
 \]
An rc-graph $R$ is \emph{associated} to a permutation $w\in\Sigma_n$ if $d(R)$ is a
reduced decomposition of $w$.
We denote the collection of rc-graphs associated to a permutation $w$ by
${\mathcal R}(w)$, and the permutation corresponding to a given rc-graph $R$
by $w(R)$.

We represent an rc-graph $R$ as a pseudoline diagram recording
the history of the inversions of $w$.
To this end, draw $n$ pseudolines going up and to the right such that the
$i$th pseudoline begins at position $(i,1)$ and ends at position $(1,w(i))$.
(The positions are numbered as in a matrix, as illustrated in
Example~\ref{exrc} below.)
The rule for constructing the pseudoline diagram is the following: two
pseudolines entering at position $(i,j)$
cross at that position if $(i,j)$ is in $R$, and otherwise
avoid each other at that position.
Note that two pseudolines cross at most once, as $d(R)$ is reduced.

Conversely, a pseudoline diagram represents a permutation $w$, where
$w(i)$ is the endpoint of the pseudoline beginning at position
$(i,1)$.
Given a pseudoline diagram where no two pseudolines cross more
than once, the crossings give an rc-graph associated to the permutation
represented by the pseudoline diagram.

\begin{exa}\label{exrc}
 Here are the pseudoline diagrams of two rc-graphs associated to the
 permutation $w=215463$.
\[
  \setlength{\unitlength}{1.pt}  
   \begin{picture}(100,100)
    \put(8,3){\epsffile{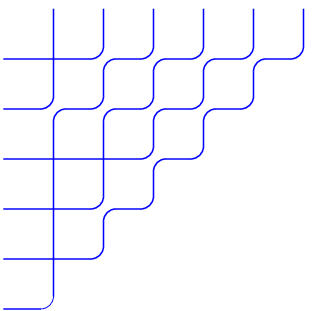}}  
    \put(0,0){6}\put(0,14.5){5}\put(0,29){4}
    \put(0,43.5){3}\put(0,58){2}\put(0,72.5){1}
    \put(20,93){1}\put(34.5,93){2}\put(49,93){3}
    \put(63.5,93){4}\put(78,93){5}\put(93.5,93){6}
   \end{picture}
   \qquad
   \begin{picture}(100,100)
    \put(8,3){\epsffile{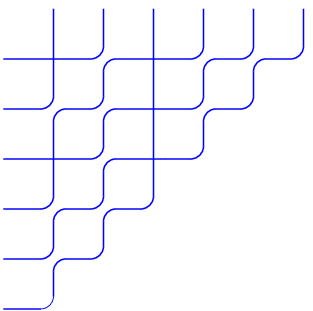}}  
    \put(0,0){6}\put(0,14.5){5}\put(0,29){4}
    \put(0,43.5){3}\put(0,58){2}\put(0,72.5){1}
    \put(20,93){1}\put(34.5,93){2}\put(49,93){3}
    \put(63.5,93){4}\put(78,93){5}\put(93.5,93){6}
   \end{picture}
 \]
\end{exa}

Given an rc-graph $R$, define the monomial
 \[
   x^R\ :=\ \prod_{(i,j)\in R} x_i\,.
\]
As shown in~\cite{BJS93,FK96},    
the Schubert polynomial $\fs_w$ indexed
by $w$ can be expressed as
 \begin{equation}\label{schub1}
   \fs_w\ =\ \sum_{R\in{\mathcal R}(w)}x^R\,.
 \end{equation}
Comparing this expression with Proposition \ref{prop:Schub-constr} suggests
that there should be a natural bijection between ${\mathcal R}(w)$ and
${\Gamma}(w,w_0)$, where $x^R\cdot x^\gamma=x^\delta$, when $R$ corresponds to
$\gamma$.
Indeed, we give such a bijection.

Given an rc-graph $R$ in ${\mathcal R}(w)$, we greedily construct the
sequence of rc-graphs $R=R_0,R_1,\ldots,R_{\ell(w_0)-\ell(w)}$ as follows.
Given $R_i$, add the pair $(k,b)$ to $R_i$ where  $(k,b)$ is the
pair not in $R_i$ that is minimal in the lexicographic order on pairs.
We have the following lemma.

\begin{lemma}\label{L:RC-chain}
 Given $R\in {\mathcal R}(w)$, construct the sequence
 $R=R_0,R_1,\ldots,R_{\ell(w_0)-\ell(w)}$ as above.
 Then
 \begin{enumerate}
  \item Every subset $R_i$ is an rc-graph.
  \item $w=w(R_0)\lessdot w(R_1)\lessdot \ldots \lessdot
          w(R_{\ell(w_0)-\ell(w)})=w_0$ is a saturated chain in Bruhat order.
  \item If $R_{i+1}$ is obtained from $R_i$ by adding the pair $(k,b)$,
       then in the labeled Bruhat order we have the labeled cover
    \[
        w(R_i)\ \xrightarrow{\ (k,b)\ }\ w(R_{i+1})\,.
    \]
  \item The labeled chain with the labels of $(3)$ is increasing.
 \end{enumerate}
\end{lemma}

\begin{exa}
 Here is an rc-graph, its pseudoline diagram, and the
 associated increasing labeled chain.
 \[
%
%
  \setlength{\unitlength}{1.pt}  
  \begin{array}{ccccc}
   \begin{array}{ccc}\mbox{ }\\\cdot&+&+\\\cdot&\cdot\\+\\\mbox{ }\end{array}
    & &
     \begin{picture}(50,30)(0,22)
      \put(0,0){\epsffile{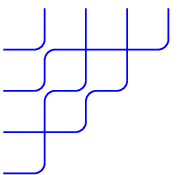}}
     \end{picture}
    & &
    1432 \xrightarrow{\,(1,1)\,} 4132 \xrightarrow{\,(2,1)\,}
    4231 \xrightarrow{\,(2,2)\,} 4321
  \end{array}
 \]
\end{exa}

\begin{proof}[Proof of Lemma~\ref{L:RC-chain}]
 We prove the first statement by induction.
 Suppose that $R_i$ is an rc-graph and we add the pair $(k,b)$ to obtain
 $R_{i+1}$, where $(k,b)$ is the minimal pair not in $R_i$.
 Consider the pseudoline diagram of $R_i$ near the position $(k,b)$:
\[
  \setlength{\unitlength}{1.pt}  
   \begin{picture}(82,50)
    \put(8,0){\epsffile{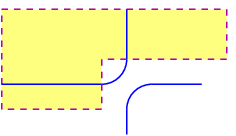}}  
    \put(0,13){$k$}        \put(43,42){$b$}
   \end{picture}
 \]
 By construction, $R_i$ contains every crossing in the shaded region, so
 the two pseudolines at this position never cross, with one
 connecting $k$ to $b$ as drawn.
 Adding the crossing gives a new pseudoline diagram for a permutation with
 exactly one more inversion. Indeed, assume that we create a
 double crossing between one of the
 two pseudolines obtained from the drawn ones by adding the extra crossing
 and another pseudoline.
 It is easy to check graphically that this can only happen if the latter
 pseudoline meets the pseudoline connecting $k$ to $b$ in $R_i$ twice. But
 this is impossible, as $R_i$ is an rc-graph.

 These arguments show that $w(R_i)\lessdot w(R_{i+1})$, implying
 the second statement.
 For the third, note that $(k,b)$ is a possible label, and the fourth
 follows by our greedy algorithm for adding pairs not in $R$.
\end{proof}

Observe that the increasing chain $\gamma(R)$ constructed in
Lemma~\ref{L:RC-chain} from an rc-graph $R$ has the property that
 \begin{equation}\label{E:star}
  \text{if $u\xrightarrow{\,(k,b)\,}v$ is a labeled cover in $\gamma(R)$,
        then $u^{-1}v=(k,l)$, for some $l>k$.}
 \end{equation}
\emph{A priori}, given that $u\xrightarrow{\,(k,b)\,}v$ is a cover, we are
only guaranteed that $u^{-1}v=(j,l)$ with $j\leq k<l$.
We prove that this property holds for any increasing chain ending in $w_0$, and
also that the map $\gamma\mapsto \gamma(R)$ is a bijection from
${\mathcal R}(w)$ to ${\Gamma}(w,w_0)$.
This will immediately imply Proposition~\ref{prop:Schub-constr}.

\begin{thm}
 Let $w\in\Sigma_n$ and $\gamma$ be an increasing chain from $w$ to $w_0$.
 Then
 \begin{enumerate}
  \item $\gamma$ has property~\eqref{E:star}.
  \item The set of pairs $(k,b)$ with $k+b\leq n$ and $(k,b)$ not a label
 of a cover in $\gamma$ is an rc-graph.
 \end{enumerate}
\end{thm}

\begin{proof}
  We prove these statements by downward induction on $\ell(w)$.
  They hold trivially when $w=w_0$.
  Consider the first two steps in $\gamma$:
 \begin{equation}\label{E:Chain}
    w\ \xrightarrow{\,(j,a)\,}\ v\ \xrightarrow{\,(k,b)\,}\ \dotsb
 \end{equation}
 Since $\gamma$ is increasing, we have $j<k$ or $j=k$ and $a<b$.
 Since the part of the chain beginning with $v$ is an increasing chain, our
 inductive hypothesis implies that subsequent covers involve only positions
 of $v$ greater than or equal to $k$.
 Thus if $i<k$, we must have $v(i)=n+1-i$ and $v(k)=b$.

 Suppose that $w^{-1}v=(i,l)$ so that $w(i)=a=v(l)$ and $i\leq j<l$.
 We consider cases $i=k$, $i=k-1$, and $i<k-1$ for $i$ separately, showing
 that in each case, we have $j=i$.

 If $i=k$, then, as $i\leq j \leq k$, we must have $j=k=i$.

 If $i<k-1$, then the betweenness condition for covers in the Bruhat order and
 the values of $v$ on $1,2,\dotsc,k{-}1$ force $l=i+1$ and so again we have
 $j=i$.

 Now suppose that $i=k-1$ and $j\neq i$ so that $j=k$ and we then have $a<b$.
 Then
 \[
   v(i)\ =\ v(k{-}1)\ =\ n{+}2{-}k\ >\ v(k)\ =\ b\ >\ a\ =\ v(l)\,,
 \]
 and so the betweenness condition on covers $w<v$ is violated.
 This completes the proof of Statement 1.

 For the second statement, consider the procedure of successively removing
 labels $(k,b)$ of covers of $\gamma$ from the rc-graph of $w_0$ (which
 contains all crossings).
 Suppose that this process applied to the chain~\eqref{E:Chain} above $v$
 has created an rc-graph $R$ associated to $v$.
 We argue that removing the crossing $(j,a)$ from $R$ creates an rc-graph
 associated to $w$.
 First, $(j,a)$ is a crossing of $R$, as we have only removed crossings that
 exceed $(j,a)$ in the lexicographic order.
 Consider the pseudoline diagram of $R$ near $(j,a)$, which we
 represent on the left.
\[
  \setlength{\unitlength}{1.pt}  
   \begin{picture}(82,50)
    \put(8,0){\epsffile{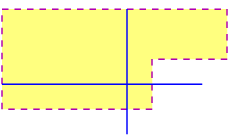}}  
    \put(0,13){$j$}        \put(43,42){$a$}
   \end{picture}
   \qquad\qquad
   \begin{picture}(82,50)
    \put(8,0){\epsffile{Rk.eps}}  
    \put(0,13){$j$}        \put(43,42){$a$}
   \end{picture}
 \]
 Note that $R$ contains all crossings in the shaded region.
 Removing the crossing at $(j,a)$ creates the picture on the right.
 Arguing as in the proof of Lemma~\ref{L:RC-chain} (but in reverse) shows
 that we obtain an rc-graph associated to $w$.
\end{proof}

\begin{rem}
 Property~\eqref{E:star} of increasing chains from $w$ to $w_0$ implies
 that the branching in the tree of increasing chains from $w$ to $w_0$ is
 quite simple; it only depends upon the permutation at a node and not on the
 history of the chain.
 The branches at a permutation $u$ consists of all possible covers
 $u\xrightarrow{(k,b)} u(k,l)$, where $k$ is the minimal position
 $i$  such that $u(i)+i<n+1$ and $b=u(i)$ and
 $l$ is any position greater than $k$ such that $u(l)>b=u(k)$ and if $k<i<l$,
 then $u(i)$ is not between $u(k)$ and $u(l)$.

 We may order the branches at $u$ by the position $l$.
 This leads to an efficient lexicographic search of this tree
 for generating all such increasing chains, and thus all rc-graphs in
 ${\mathcal R}(w)$ as well as the multiset of monomials in the
 Schubert polynomial $\frakS_{w}(x)$.
 The computational complexity and memory usage of this algorithm
 depend only upon $n$ (the number of entries in $w$), the length $l(w)$ of
 $w$, and the number of increasing chains, $c:=\frakS_{w}(1,1,\dotsc,1)$.
 This is true, in fact, for any  algorithm generating $c$ combinatorial
 objects corresponding to the monomials in the Schubert polynomial
 $\frakS_{w}(x)$. Furthermore, it is not hard to see that the above
 algorithm is of order $O(nlc)$, where $l:=l(w_0)-l(w)=\binom{n}{2}-l(w)$.
 Indeed, note that the number of internal nodes of a rooted tree with a
 fixed number of leaves ($c$ in our case) and a fixed distance from its
 root to all its leaves ($l$ in our case) is maximized by the tree with all
 internal nodes having only one descendant; thus, there are $O(lc)$ nodes.
 On the other hand, finding the descendants of a given node is of order
 $O(n)$.

 Let us also note that a closely related algorithm, with the same complexity
 as the one above, appears in~\cite{KM01,Mil02}; it is described in terms of
 removing crossings from the rc-graph corresponding to $w_0$.
 More precisely, the authors of the mentioned papers use combinatorial
 versions of the divided difference operators in order to recursively
 generate the rc-graphs corresponding to a permutation $w$.
 One starts from $w_0$ and applies these operators using any chain in
 the weak order on $\Sigma_n$ from $w_0$ to $w$.
 A certain choice of chain guarantees that each rc-graph in the obtained tree
 has at least one descendant upon applying the corresponding operator.
 It turns out that the rc-graphs on a given level $i$ in this tree
 correspond (via the bijection  described above) to the subchains starting
 at level $l-i$ in the tree constructed by the previous algorithm.

 While this algorithm for generating the multiset of monomials in the
 Schubert polynomial $\mathfrak{S}_w(x)$ also computes
 this polynomial, A.~Buch has pointed out to us that the transition
 formula of Lascoux and Sch\"utzenberger~\cite{LS82a} provides a more
 efficient algorithm to {\it compute} $\mathfrak{S}_w(x)$.
\end{rem}

We conclude by mentioning the connection with the insertion
algorithm for rc-graphs in \cite{KK00}, which bijectively proves the Pieri
formula~\eqref{E:Pieri} for multiplying a Schubert
polynomial by a complete homogeneous symmetric polynomial
$h_i(x_1,\ldots,x_k)$.
Given an rc-graph $R$ with $x^R=x_1^{\alpha_1}\ldots
x_{n-1}^{\alpha_{n-1}}$, the chain $\gamma(R)$ can be obtained by successive
insertions into $R$ corresponding to multiplications of $x^R$ by the monomials
$x_k^{n-k-\alpha_k}$ in $h_{n-k-\alpha_k}(x_1,\ldots,x_k)$, for
$k=1,\ldots,n-1$.
In this this special case, the insertion algorithm
reduces to the older one in \cite{BeBi}.\smallskip

We thank Ezra Miller for his detailed comments on an earlier version of this
paper.

\def\cprime{$'$}
\providecommand{\bysame}{\leavevmode\hbox to3em{\hrulefill}\thinspace}
\providecommand{\MR}{\relax\ifhmode\unskip\space\fi MR }
\providecommand{\MRhref}[2]{%
  \href{http://www.ams.org/mathscinet-getitem?mr=#1}{#2}
}
\providecommand{\href}[2]{#2}

\end{document}